\documentclass[journal]{IEEEtran}
\usepackage{amsmath,amssymb,amsfonts}
\usepackage{cases}
\usepackage{color,xcolor}
\usepackage{graphicx}
\usepackage{subfigure}
\usepackage{epstopdf}
\usepackage{multirow}
\usepackage{rotating}
\usepackage{booktabs}
\usepackage{tikz}
\usepackage{acronym}
\usepackage{epstopdf}
\usepackage{pifont} % 导入pifont包以使用打钩和打叉符号
\usepackage{makecell} % 表格换行
\usepackage[most]{tcolorbox}
\definecolor{boxdarkgray}{RGB}{85,85,85}
\newcommand{\best}[1]{\textcolor{red}{#1}}
\newcommand{\second}[1]{\textcolor{blue}{#1}}
\usepackage{amssymb}
\usepackage{tabularx}
\usepackage[table]{xcolor}
\usepackage[ruled,vlined,linesnumbered]{algorithm2e}
\SetKwInput{KwInput}{Input}
\SetKwInput{KwOutput}{Output}

\definecolor{OptBlue}{RGB}{0,0,150}
\definecolor{OptOrange}{RGB}{190,65,0}
\definecolor{OptPurple}{RGB}{120,0,160}
\definecolor{OptGreen}{RGB}{0,100,0}

\usepackage{xcolor}

\definecolor{OptBlue}{RGB}{0,20,105}
\definecolor{OptOrange}{RGB}{185,70,0}
\definecolor{OptPurple}{RGB}{115,0,145}
\definecolor{OptGreen}{RGB}{0,85,25}

\definecolor{OptBlueBg}{RGB}{226,234,245}
\definecolor{OptOrangeBg}{RGB}{250,238,218}
\definecolor{OptPurpleBg}{RGB}{240,230,245}
\definecolor{OptGreenBg}{RGB}{229,242,224}

\newtcolorbox{mybox}[1][]{
    colback=white,
    colframe=gray!80,
    colbacktitle=gray!80,
    coltitle=white,
    fonttitle=\bfseries\itshape\large,
    title={\centering #1},
    boxrule=0.8pt,
    titlerule=0pt,
    arc=2pt,
    left=5pt,
    right=5pt,
    top=5pt,
    bottom=5pt,
    toptitle=3pt,
    bottomtitle=3pt
}

\usepackage[demo,
            export]{adjustbox}  % it also load graphicx
\usepackage{tabularray}
\usepackage{colortbl}
\usepackage[colorlinks,
linkcolor=blue,
anchorcolor=blue,
citecolor=blue
]{hyperref}
\usepackage{amssymb}
\usepackage{array} % 引入array宏包
\usepackage{catchfilebetweentags}
\definecolor{tableTitleColor}{rgb}{0.8,0.8,0.8}
\begin{document}

\title{OptGraph: Large Language Models Enhanced Evolutionary Optimization Via Graph Retrieval-Augmented Generation}

\author{Xianchao Xiu, Jianhao Li, Huangyue Chen, and Wanquan Liu, \IEEEmembership{Senior Member,~IEEE}

\thanks{This work was supported in part by the National Natural Science Foundation of China (12371306, 12501452). (\textit{Corresponding author: Huangyue Chen}.)}
\thanks{Xianchao Xiu and Jianhao Li are with the School of Mechatronic Engineering and Automation, Shanghai University, Shanghai 200444, China (e-mail: xcxiu@shu.edu.cn; lijianhao@shu.edu.cn).}
\thanks{Huangyue Chen is with the School of Mathematics \& Center for Applied Mathematics of Guangxi, Guangxi University, Nanning, 530004, China (e-mail: hychen2024@gxu.edu.cn).}
\thanks{Wanquan Liu is with the School of Intelligent Systems Engineering, Sun Yat-sen University, Guangzhou 510275, China (e-mail: liuwq63@mail.sysu.edu.cn).}
}

\maketitle

\begin{abstract}
Large language models (LLMs) have emerged as a powerful tool for automated evolutionary optimization, but existing methods remain limited in pattern reuse,  error-aware refinement, and retrieval robustness across diverse tasks. To address these limitations, we propose OptGraph, the first optimization agentic workflow that introduces graph retrieval-augmented generation (GraphRAG). Specifically, OptGraph first constructs reusable experience as a typed graph, capturing the relationships among modeling patterns, problem formalization, implementation details, and error corrections. In the inference stage, OptGraph leverages graph neighborhood information to enrich retrieved knowledge, providing structured context to improve modeling, verification, and iterative refinement. Moreover, OptGraph supports adaptive knowledge updates, enabling the distillation of execution traces and verification feedback into reusable graph knowledge without ndertaking LLM parameter tuning. Extensive experiments on benchmark datasets show that our proposed OptGraph achieves an average exact accuracy 8.9\% higher than the state-of-the-art prompt-based automated optimization frameworks. 
Our code has been made available at \href{https://github.com/xianchaoxiu/OptGraph}{https://github.com/xianchaoxiu/OptGraph}.
\end{abstract}


\section{Introduction}

\IEEEPARstart{O}{ptimization} problems are fundamental to modern decision-making systems, with broad applications spanning transportation, finance and manufacturing, see the recent surveys \cite{archetti2022optimization,gunjan2023brief,weckenborg2024flexibility}.
In practice, the bottleneck often lies not only in solving the optimization model but also in formulating it. Practitioners must translate ambiguous natural-language descriptions into decision variables, objectives, and constraints, and subsequently implement the resulting model as executable code compatible with optimization solvers such as Gurobi \cite{gurobi2026}  and CPLEX \cite{ibm2026cplex}. Clearly, this process needs extensive expertise in operations research (OR), which makes optimization both time-consuming and difficult to directly apply to diverse real-world scenarios \cite{wu2024evolutionary,tiande2025machine,fan2025artificial}.

During the last few years, large language models (LLMs) based on the Transformer architecture \cite{vaswani2017attention} have shown great potential in evolutionary optimization and automated optimization \cite{ramamonjison2023NL4Opt,chen2024diagnosing,chen2026solver}. With emerging capabilities in mathematical reasoning \cite{ahn2024large},
program synthesis \cite{ma2026llamoco}, and algorithm design \cite{liu2024llm4ad}, LLMs are able to interpret problem semantics, construct logically coherent mathematical formulations, and generate solver-ready implementations. 
According to  \cite{zhang2025systematic,xiu2026large}, existing methods primarily fall into two categories: learning-based and prompt-based. Broadly speaking, learning-based methods adapt LLMs using optimization-specific data, execution feedback, or reinforcement learning, and representatives include  ORLM \cite{huang2025orlm}, OptMATH \cite{lu2025optmath}, LLMOPT \cite{jiang2025llmopt}, and OR-R1 \cite{ding2026or}. In contrast, prompt-based methods construct reasoning workflows, retrieval strategies, agent collaboration, or validation mechanisms around pretrained LLMs, and typical examples are  Chain-of-Experts (CoE) \cite{xiao2024chain}, OptiMUS \cite{ahmaditeshnizi2024optimus0.2}, OptiTree \cite{liu2025optitree}, and Lean-LLM-OPT \cite{liang2026large}. All these methods reflect a broader evolution in optimization research, extending beyond classical  iterative optimization toward LLM-powered automated optimization \cite{zheng2026survey}. By significantly reducing the runtime and human effort required for optimization modeling and algorithm design, these methods lower the barrier to using optimization techniques \cite{shi2026mirror,huang2025autonomous,ge2025mora}.

\begin{figure*}[!t]
    \centering
    \subfigure[]{
        \begin{minipage}{0.33\textwidth}
            \centering
            \includegraphics[width=\linewidth]{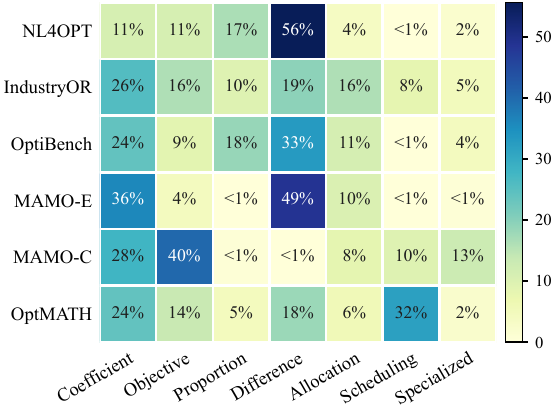}
            \vspace{-2mm} 
            \label{fig:pattern_coverage}
        \end{minipage}
    }
    \hfill
    \subfigure[]{
        \begin{minipage}{0.28\textwidth}
            \centering
            \includegraphics[width=\linewidth]{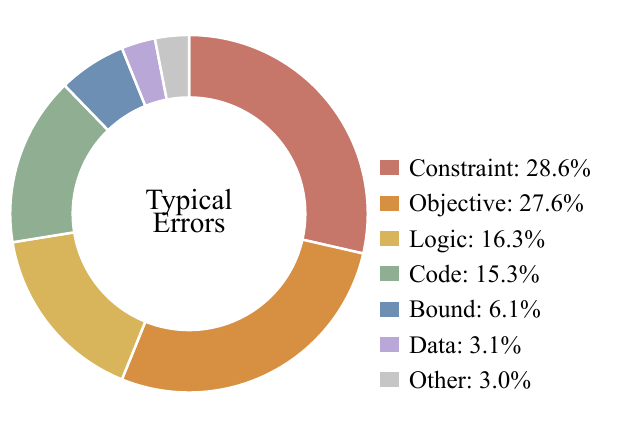}
           \vspace{4mm}
            \label{fig:typical_errors}
        \end{minipage}
    }
    \hfill
    \subfigure[]{
        \begin{minipage}{0.33\textwidth}
            \centering
            \includegraphics[width=\linewidth]{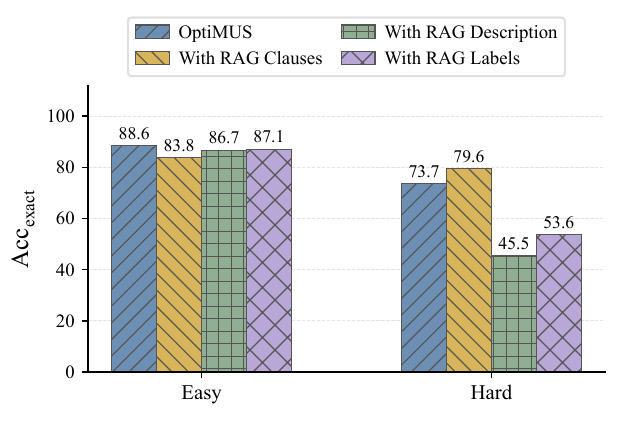}
            \vspace{-2mm}
            \label{fig:vector_rag_comparison}
        \end{minipage}
    }

   % \vspace{-2mm}
    \caption{Observations for existing automated optimization frameworks, where  (a) diverse reusable patterns across datasets, (b) recurring formulation errors, (c) unstable effects of naive RAG for OptiMUS.}
    \label{fig:motivating_observations}
\end{figure*}

Even though LLM-driven optimization has achieved tremendous success, it still faces three drawbacks. To begin with,  it is shown in Fig. \ref{fig:pattern_coverage} that optimization benchmark datasets involve heterogeneous yet reusable modeling patterns, such as objective-related, proportional, difference, allocation, scheduling, and specialized OR structures. This suggests that reusable modeling knowledge should be organized around problem types and variants, rather than stored as isolated examples or static templates \cite{kong2025alphaopt}. Furthermore, as illustrated in Fig. \ref{fig:typical_errors}, LLM-based failures are not limited to code execution errors, but also contain recurring formulation mistakes, such as objective errors, missing constraints, wrong constraint logic, bound errors, and data-interpretation errors. This highlights the need for a validation mechanism that checks objective consistency, constraint completeness, and program correctness, together with error-aware knowledge for repair \cite{lian2026reloop}. Last but not least, Fig. \ref{fig:vector_rag_comparison} reports the naive retrieval-augmented generation (RAG) \cite{lewis2020retrieval} results from OptiMUS \cite{ahmaditeshnizi2024optimus0.3}, where retrieved examples can help hard problems but may also degrade performance when the retrieved context is less relevant. This demonstrates that retrieval robutness is critical for optimization, since textual similarity does not always match modeling logic, which has been explored in \cite{zhu2026llm}. Overall, automated   optimization requires structured knowledge for pattern reuse, validation-based error correction, and robust retrieval across diverse optimization tasks.

In fact, the aforementioned issues are largely attributable to hallucination prevalent in LLMs when dealing with mathematical reasoning, domain knowledge, and executable code \cite{huang2025survey}. Popular methods such as Chain-of-Thought \cite{wei2022chain}, Tree-of-Thought \cite{yao2023tree}, and Graph-of-Thought \cite{besta2024graph} improve reliability through structured reasoning. In addition, RAG incorporates external knowledge, but naive RAG mainly retrieves isolated text chunks based on semantic similarity \cite{chen2024benchmarking}. Graph-augmented retrieval instead organizes knowledge and its relationships in graph structures. Graph RAG (GraphRAG) \cite{edge2024local,peng2025graph}, together with  HippoRAG \cite{gutierrez2024hipporag} and LightRAG \cite{guo2024lightrag}, validates the potential of this paradigm. \textit{This raises a natural question: can graph-augmented retrieval boost automated optimization modeling by connecting problem patterns, formulations, code snippets, typical errors, and repair hints?}

To fill this gap, we develop an adaptive GraphRAG-based agentic workflow for automated optimization, dubbed OptGraph. To the best of our knowledge, it is the first framework inspired by graph-augmented retrieval designed for optimization modeling and solving. Compared to existing prompt-based frameworks such as OptiTree \cite{liu2025optitree} and Lean-LLM-OPT \cite{liang2026large}, it achieves superior accuracy and robustness. 
The main contributions of this paper are as follows:
\begin{itemize}
    \item We characterize reusable optimization knowledge as a heterogeneous graph and retrieve structurally connected knowledge to facilitate formulation and  code generation within the framework.
    
    \item We introduce a validation-guided repair mechanism that evaluates objective consistency, constraint completeness, and program correctness, and leverages error-aware graph knowledge to correct failed formulations and implementations through iterative refinement.
    
    \item We design an adaptive graph update that distills execution traces and validation feedback into reusable knowledge, enabling non-parametric accumulation across tasks.
\end{itemize}

The remainder of this paper is organized as follows.
Section \ref{sec:related_work} reviews LLM-driven optimization frameworks.
Section \ref{sec:methodology} details our OptGraph, covering GraphRAG construction, knowledge modeling, validation mechanism, and GraphRAG update.
Section \ref{sec:experiments} provides experimental results and ablation studies.
Finally, Section \ref{sec:conclusion} concludes this paper with some research directions.

\begin{figure*}[!t]
\centering
\includegraphics[width=\textwidth]{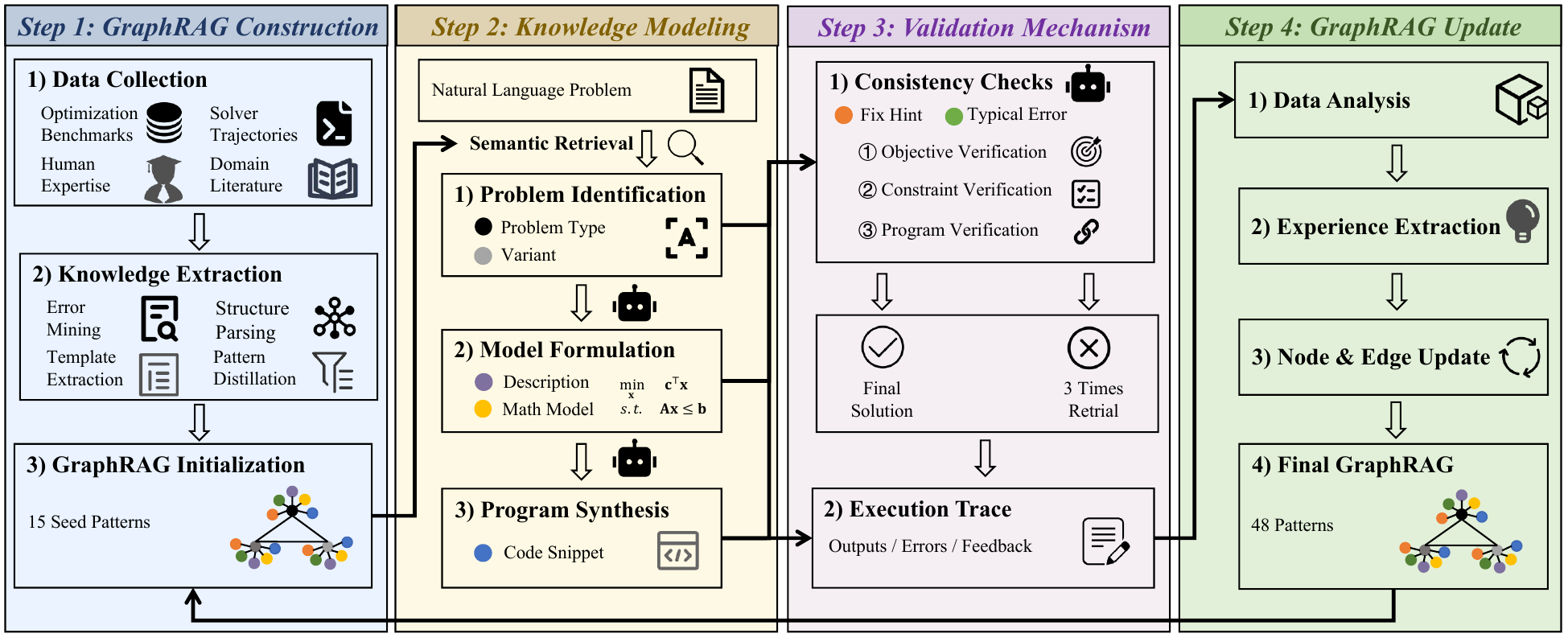}
\caption{Overall workflow of the proposed OptGraph. Given a natural-language optimization problem, OptGraph retrieves structured knowledge from a heterogeneous modeling graph and guides LLM agents to generate a mathematical formulation and executable code. The generated code is then executed and checked by the validation mechanism. If the result is not accepted, OptGraph retrieves error-aware knowledge and performs repair.}
\label{fig:optgraph_workflow}
\end{figure*}

\section{Related Work}
\label{sec:related_work}

This section reviews some LLM-based automated optimization and discusses their differences from our OptGraph.

\subsection{Learning-Based Methods}

As previously mentioned, learning-based methods leverage optimization-specific data to enable LLMs for automated optimization. ORLM \cite{huang2025orlm} trains open-source LLMs on OR-Instruct and introduces a new dataset called  IndustryOR. OptMATH \cite{lu2025optmath} constructs optimization modeling data through bidirectional synthesis and filtering.
Additionally, LLMOPT \cite{jiang2025llmopt} employs instruction tuning, model alignment, and self-correction. More recently, OR-R1 \cite{ding2026or} combines supervised fine-tuning with execution-based reinforcement learning, while AutoOR \cite{motwani2026autoor} develops a synthetic data generation and reinforcement learning pipeline to automate the formalization of linear, mixed-integer, and nonlinear optimization problems. Unlike these methods, our OptGraph organizes external knowledge and experience in an adaptive updated graph without modifying the underlying LLM parameters.

\subsection{Prompt-Based Methods}

Prompt engineering has been widely used in LLMs, including automated optimization. For example, CoE \cite{xiao2024chain} introduces a multi-agent workflow that combines forward construction with backward reflection. OptiMUS \cite{ahmaditeshnizi2024optimus0.2,ahmaditeshnizi2024optimus0.3} decomposes long problem descriptions into objectives and constraints and incorporates retrieval-augmented generation. Recently, OptiTree \cite{liu2025optitree} organizes optimization knowledge into hierarchical problem categories, while MCTS-based Autoformulation \cite{astorga2025autoformulation} searches over candidate variables, constraints, and objectives using value estimation and symbolic pruning. Lean-LLM-OPT \cite{liang2026large} combines few-shot prompting, retrieval, and agent collaboration to address complex optimization problems. Other studies introduce validation, debugging, and reusable skills to improve formulation reliability \cite{yang2026optskills,liu2026opt}. 
These methods validate the potential of pretrained LLMs for automated optimization modeling. However, most still rely on predefined reasoning procedures, static knowledge structures, or isolated retrieved examples, rather than adaptively accumulating and reusing structured modeling knowledge across tasks.

\section{Methodology}
\label{sec:methodology}

This section presents the architecture and implementation details of our proposed OptGraph. Specifically, it contains four main components: GraphRAG construction, knowledge modeling, validation mechanism, and GraphRAG update. 
The overall workflow is illustrated in Fig. \ref{fig:optgraph_workflow}.

\subsection{GraphRAG Construction}

In Step 1, OptGraph constructs an initial graph-structured knowledge base from reusable optimization resources.

\subsubsection{Data Collection}
OptGraph first collects prior knowledge from optimization benchmarks, solver trajectories, human expertise, and domain literature. These sources provide natural-language problem descriptions, reference formulations, executable code, execution outputs, common modeling errors, and repair experience.

\subsubsection{Knowledge Extraction}

Through error mining, structure parsing, template extraction, and pattern distillation, OptGraph transforms these raw modeling resources into graph-structured knowledge. Specifically, error mining identifies recurring formulation and program errors, structure parsing organizes records into typed knowledge units, template extraction summarizes reusable formulation and code skeletons, and pattern distillation merges similar cases into compact modeling patterns. Mathematically, the extracted knowledge is represented as the following  directed heterogeneous graph
\begin{equation}
\mathcal{G} = (\mathcal{V}, \mathcal{E}, \tau, \rho),
\label{eq:graph_def}
\end{equation}
where $\mathcal{V}$ and $\mathcal{E}$ denote the sets of nodes and edges, respectively, and $\tau$ and $\rho$ map nodes and edges to their corresponding types.
The node types include
\begin{equation}
\mathcal{T}_v =
\left\{
\begin{array}{l}
\texttt{Problem Type},\ \texttt{Variant},\\
\texttt{Description},\ \texttt{Math Model},\\
\texttt{Code Snippet},\  \texttt{Typical Error},  \\
\texttt{Fix Hint}
\end{array}
\right\}.
\label{eq:node_types}
\end{equation}

\subsubsection{GraphRAG Initialization}
The output of Step 1 is an initial GraphRAG containing seed modeling patterns that connect the \texttt{Problem Type} with \texttt{Variant}, \texttt{Description}, \texttt{Math Model}, \texttt{Code Snippet}, \texttt{Typical Error}, and \texttt{Fix Hint}. 
To make the graph structure concrete, Fig. \ref{fig:modeling_graph_example} gives an illustrative example. It is seen that such a structure preserves the relations among problem semantics, formulations, implementations, and repair knowledge, enabling structured retrieval before graph updates are performed.

\begin{figure*}[!t]
\centering
\small
\begin{mybox}[\small An Example of the Modeling Graph]
\label{box:modeling_graph_example}
\textbf{\texttt{Problem Type:}} Vehicle routing problem (VRP)

\textbf{\texttt{Variant:}} Capacitated vehicle routing problem (CVRP)

\textbf{\texttt{Description:}} Design vehicle delivery routes to serve all customers with minimal travel cost under vehicle capacity constraints

\textbf{\texttt{Math Model:}}
\[
\begin{aligned}
\min_{\mathbf{x}} \quad
& \sum_{k \in K}\sum_{i \in V}\sum_{j \in V}
c_{ij}x_{ijk} \\
\text{s.t.} \quad
& \sum_{k \in K}\sum_{j \in V}x_{ijk}=1,
&& \forall i \in C \\
& \sum_{i \in C} d_i
   \sum_{j \in V}x_{ijk}\leq Q,
&& \forall k \in K \\
& x_{ijk}\in\{0,1\},
&& \forall i,j\in V,\ k\in K
\end{aligned}
\]

\textbf{\texttt{Code Snippet:}}
\[
\begin{array}{l}
\texttt{x = addVars(V, V, K, vtype=GRB.BINARY)}\\
\texttt{min sum(c[i,j] * x[i,j,k] for i in V for j in V for k in K)}\\
\texttt{addConstrs(each customer visited once), addConstrs(vehicle capacity limits)}
\end{array}
\]

\textbf{\texttt{Typical Error:}}
Missing customer coverage, ignoring vehicle capacity, disconnected subtours

\textbf{\texttt{Fix Hint:}} Visit each customer once, add capacity constraints, eliminate subtours
\end{mybox}
\caption{A simple example of the modeling graph in GraphRAG initialization, where all displayed entries are treated as nodes connected by semantic and functional edges.}
\label{fig:modeling_graph_example}
\end{figure*}

\subsection{Knowledge Modeling}
In Step 2, OptGraph utilizes the initial graph-structured knowledge base to guide the modeling process. Let $\mathcal{V}_{\mathrm{var}}$ denote the set of all \texttt{Variant} nodes and each node represents a reusable modeling pattern. 
Given a natural-language optimization problem $\mathcal{P}$, OptGraph computes a relevance score between the problem and each variant node by
\begin{equation}
    S(\mathcal{P}, v),
    \quad v \in \mathcal{V}_{\mathrm{var}},
\label{eq:retrieval_score}
\end{equation}
where $S(\mathcal{P},v)$ measures the semantic and modeling-pattern relevance between the input problem and node $v$. 
The top-$k$ \texttt{Variant}  nodes are selected via
\begin{equation}
    \mathcal{R}_k(\mathcal{P})
    =
    \operatorname{TopK}_{v \in \mathcal{V}_{\mathrm{var}}}
    S(\mathcal{P}, v).
\label{eq:topk_retrieval}
\end{equation}
Subsequently, OptGraph expands the graph neighborhood of the retrieved \texttt{Variant} nodes as follows
\begin{equation}
    \mathcal{K}(\mathcal{P})
    =
    \bigcup_{v \in \mathcal{R}_k(\mathcal{P})}
     (\{v\} \cup \mathcal{N}(v) ),
\label{eq:retrieved_context}
\end{equation}
where $\mathcal{N}(v)$ represents the adjacent nodes of $v$. 
The retrieved context $\mathcal{K}(\mathcal{P})$ may contain the related \texttt{Problem Type}, \texttt{Description}, \texttt{Math Model}, \texttt{Code Snippet}, \texttt{Typical Error}, and \texttt{Fix Hint}. 
Compared with naive RAG or conventional example retrieval, this graph expansion provides connected modeling knowledge rather than isolated text snippets.

\subsubsection{Problem Identification}
According to $\mathcal{K}(\mathcal{P})$ in Eq. \eqref{eq:retrieved_context}, the problem-identification agent determines the likely \texttt{Problem Type} and \texttt{Variant}. This process helps OptGraph select appropriate modeling patterns before generating the model formulation.

\subsubsection{Model Formulation}
The model-formulation agent (denoted as $A_{\mathrm{model}}$) generates a mathematical formulation by combining the input problem with the retrieved graph context. 
Specifically,
\begin{equation}
    \hat{\mathcal{M}}
    =
    A_{\mathrm{model}}(
    \mathcal{P}, \mathcal{K}(\mathcal{P})
    ).
\label{eq:model_agent}
\end{equation}
The retrieved \texttt{Description} and \texttt{Math Model} nodes provide formulation guidance, covering decision variables, objective functions, and constraints.

\subsubsection{Program Synthesis}
Finally, the program-synthesis agent  (denoted as $A_{\mathrm{code}}$)  generates executable code, i.e.,
\begin{equation}
    \hat{\mathcal{C}}
    =
    A_{\mathrm{code}}(
    \mathcal{P}, \hat{\mathcal{M}}, \mathcal{K}(\mathcal{P})
    ).
\label{eq:code_agent}
\end{equation}
The retrieved \texttt{Code Snippet} nodes provide implementation skeletons, while related \texttt{Typical Error} and \texttt{Fix Hint} nodes mitigate common formulation-code mismatches. 
The generated formulation and code are then passed to the validation mechanism.

\begin{algorithm*}[t]
\DontPrintSemicolon
\caption{Workflow of Our Proposed OptGraph}
\label{alg:optgraph_modeling}

\KwInput{Optimization problem set
$\{\mathcal{P}_q\}_{q=1}^{Q}$, initial graph $\mathcal{G}^{(0)}$,
maximum repair rounds $T$, retrieval size $k$}
\KwOutput{Final answers $\{\hat{y}_q\}_{q=1}^{Q}$,
execution traces $\{\xi_q\}_{q=1}^{Q}$,
updated graph $\mathcal{G}^{(Q)}$}

\BlankLine
\textcolor{OptBlue}{
\textbf{\texttt{// Step 1: GraphRAG Construction}}
}\;
Initialize GraphRAG with the heterogeneous graph
$\mathcal{G}^{(0)}$ defined in Eq. \eqref{eq:graph_def}\;

\For{$q=1$ \KwTo $Q$}{

    \BlankLine
    \textcolor{OptOrange}{
    \textbf{\texttt{// Step 2: Knowledge Modeling}}
    }\;
    
    Retrieve the top-$k$ variant nodes from
    $\mathcal{G}^{(q-1)}$ according to
    Eq. \eqref{eq:topk_retrieval}\;
    
    Expand their graph neighborhoods and obtain the retrieved context
    $\mathcal{K}(\mathcal{P}_q)$ according to
    Eq. \eqref{eq:retrieved_context}\;
    
    Identify the problem type and relevant modeling variants from
    $\mathcal{K}(\mathcal{P}_q)$\;
    
    Generate the initial formulation
    $\hat{\mathcal{M}}_q^{(0)}$ and executable code
    $\hat{\mathcal{C}}_q^{(0)}$ according to
    Eqs. \eqref{eq:model_agent} and \eqref{eq:code_agent}\;

    \BlankLine
    \textcolor{OptPurple}{
    \textbf{\texttt{// Step 3: Validation Mechanism}}
    }\;
    
    Set $r=0$, $\hat{y}_q=\varnothing$, and
    $\xi_q=\varnothing$\;

    \While{$r \le T$}{

        Execute $\hat{\mathcal{C}}_q^{(r)}$ and obtain
        $(o_q^{(r)},\hat{y}_q^{(r)},s_q^{(r)})$
        according to Eq. \eqref{eq:execution}\;
        
        Define the current modeling and execution state
        $\eta_q^{(r)}$ according to
        Eq. \eqref{eq:repair_state}\;
        
        Validate the formulation, code, and answer to obtain
        $z_q^{(r)}$ and validation feedback $f_q^{(r)}$
        according to Eq. \eqref{eq:validation}\;
        
        Append
        $(\eta_q^{(r)},z_q^{(r)},f_q^{(r)})$
        to $\xi_q$ according to
        Eq. \eqref{eq:execution_trace}\;

        \If{$z_q^{(r)}=1$}{
            Set $\hat{y}_q=\hat{y}_q^{(r)}$\;
            \textbf{break}\;
        }

        \If{$r=T$}{
            \textbf{break}\;
        }

        Retrieve error-aware context
        $\mathcal{K}_{\mathrm{err}}(\mathcal{P}_q)$
        from $\mathcal{G}^{(q-1)}$\;
        
        Repair the formulation and code according to
        Eq. \eqref{eq:repair}\;
        
        Set $r=r+1$\;
    }

    \If{$\hat{y}_q=\varnothing$}{
        Mark $\xi_q$ as a failed execution trace\;
    }

    \BlankLine
    \textcolor{OptGreen}{
    \textbf{\texttt{// Step 4: GraphRAG Update}}
    }\;
    
    Analyze the execution trace $\xi_q$\;
    
    Extract the candidate graph update
    $\Delta\mathcal{G}_q$
    according to Eq. \eqref{eq:candidate_update}\;
    
    Update the graph as
    $\mathcal{G}^{(q)}
    =
    \mathcal{G}^{(q-1)}
    \oplus
    \Delta\mathcal{G}_q$
    according to Eq. \eqref{eq:graph_update}\;
}

\BlankLine
\Return{$\{\hat{y}_q\}_{q=1}^{Q}$,
$\{\xi_q\}_{q=1}^{Q}$, and $\mathcal{G}^{(Q)}$}
\end{algorithm*}

\subsection{Validation Mechanism}

In Step 3, OptGraph executes the generated code, validates the generated result, and records the modeling and execution process for subsequent graph updating.

\subsubsection{Consistency Checks}

The generated code $\hat{\mathcal{C}}$ is first executed to obtain
\begin{equation}
    (o, \hat{y}, s)
    =
    \operatorname{Exec}(\hat{\mathcal{C}}),
\label{eq:execution}
\end{equation}
where $o$ denotes the raw execution output, $\hat{y}$ is the extracted numerical answer, and $s$ represents the execution status.

Note that the validation mechanism performs consistency checks in three dimensions: objective verification, constraint verification, and program verification. 
These checks determine whether the generated objective is consistent with the problem intent, whether the constraints cover the required conditions, and whether the generated program can be correctly executed and parsed. 
The validation result is defined as
\begin{equation}
    z
    =
    V(\mathcal{P}, \hat{\mathcal{M}},
    \hat{\mathcal{C}}, o, \hat{y}, s),
\label{eq:validation}
\end{equation}
where $z \in \{0,1\}$ indicates whether the generated result is accepted.
At the $t$-th repair round, the corresponding validation result is denoted by $z^{(t)}$.
If validation fails, i.e., $z^{(t)}=0$, OptGraph retrieves error-aware context from the graph, including relevant \texttt{Typical Error} and \texttt{Fix Hint} nodes, and invokes the repair agent.
For compact notation, the modeling and execution state at the $t$-th repair round is defined as
\begin{equation}
    \eta^{(t)}
    =
    (
    \hat{\mathcal{M}}^{(t)},
    \hat{\mathcal{C}}^{(t)},
    o^{(t)},
    \hat{y}^{(t)},
    s^{(t)}
   ).
\label{eq:repair_state}
\end{equation}

Next, the repair agent updates the formulation and code by conditioning on the original problem, the current state, and the retrieved error-aware context, given by
\begin{equation}
\begin{aligned}
    (
    \hat{\mathcal{M}}^{(t+1)},
    \hat{\mathcal{C}}^{(t+1)}
    )
    =
    A_{\mathrm{repair}}
   (
    \mathcal{P},
    \eta^{(t)},
    \mathcal{K}_{\mathrm{err}}(\mathcal{P})
    ).
\end{aligned}
\label{eq:repair}
\end{equation}
The validation-and-repair process continues until the generated result is accepted or the maximum number of repair rounds is reached. In our implementation, this maximum number is set to $3$.

\subsubsection{Execution Trace}

For each problem, OptGraph records the intermediate formulations, generated programs, execution outputs, numerical answers, execution statuses, and validation feedback across all rounds. The resulting execution trace is defined as
\begin{equation}
    \xi
    =
    \{
    (
    \eta^{(t)},
    z^{(t)},
    f^{(t)}
    )
    \}_{t=0}^{\tau},
\label{eq:execution_trace}
\end{equation}
where $f^{(t)}$ denotes the validation feedback at the $t$-th  round, and $\tau$ denotes the terminal round at which the result is accepted or the repair limit is reached.

If $z^{(\tau)}=1$, OptGraph returns the accepted result as the final solution. Otherwise, it records a failed execution trace. All execution outputs, errors, and validation feedback contained in $\xi$ are subsequently considered for experience extraction and  graph updating.

\subsection{GraphRAG Update}

In Step 4, OptGraph updates the graph after validation and repair to incorporate newly observed modeling experience. 

\subsubsection{Data Analysis}

For each problem $\mathcal{P}_q$ with $q=1,\cdots,Q$, OptGraph analyzes the execution trace $\xi_q$ defined in Eq. \eqref{eq:execution_trace}. The trace records the retrieved context, generated formulation and code, execution outputs, numerical answer, execution status, and validation feedback. Successful, failed, and repaired cases are examined to identify reusable modeling knowledge and recurring formulation errors.

\subsubsection{Experience Extraction}

OptGraph extracts reusable experience from the execution trace. The candidate graph update obtained from $\xi_q$ is defined as
\begin{equation}
    \Delta \mathcal{G}_q
    =
    U(\xi_q)
    =
    (
    \Delta \mathcal{V}_q,
    \Delta \mathcal{E}_q
   ),
\label{eq:candidate_update}
\end{equation}
where $\Delta \mathcal{V}_q$ and $\Delta \mathcal{E}_q$ denote the candidate nodes and edges, respectively. The extracted knowledge may include new modeling \texttt{Variant}, \texttt{Math Model}, \texttt{Code Snippet}, \texttt{Typical Error} and \texttt{Fix Hint} nodes, and their semantic or functional relations.

\begin{table*}[t]
\centering
\caption{Accuracy comparisons (exact accuracy/5\% tolerance accuracy) of representative methods on selected datasets, where the best and second-best results are marked in red and blue, respectively. }
\label{tab3}
\renewcommand{\arraystretch}{1.2}
\footnotesize
\setlength{\tabcolsep}{3pt}
\begin{tabularx}{\linewidth}{l *{6}{>{\centering\arraybackslash}X} >{\centering\arraybackslash}X}
\toprule
Methods
& NL4Opt
& IndustryOR 
& OptiBench 
& MAMO-E
& MAMO-C
& OptMATH
& Average \\
\midrule
GPT-5.1  (2025)
& 72.9\%/78.0\% 
& 43.0\%/51.0\% 
& 58.0\%/65.8\% 
& 81.6\%/85.8\% 
& 47.8\%/58.1\% 
& 27.7\%/30.7\%
& 55.2\%/61.6\% \\

Gemini-3.1-Pro  (2026)
& 73.4\%/78.5\%
& 42.0\%/50.0\%
& 58.2\%/66.0\%
& 82.1\%/86.1\%
& 46.8\%/57.1\%
& 27.1\%/30.1\%
& 54.9\%/61.3\% \\

DeepSeek-V3.2  (2025)
& 72.4\%/77.6\%
& 44.0\%/50.0\%
& 57.4\%/65.1\%
& 81.0\%/85.0\%
& 48.3\%/57.6\%
& 26.5\%/30.1\%
& 54.9\%/60.9\% \\

Qwen3.6-Plus  (2026)
& 69.6\%/74.8\%
& 42.0\%/49.0\%
& 55.4\%/62.6\%
& 78.0\%/81.8\%
& 46.3\%/55.7\%
& 27.1\%/28.9\%
& 53.1\%/58.8\% \\

\midrule
CoE   (2024)
& 77.1\%/83.6\% 
& 51.0\%/61.0\% 
& \second{65.5\%}/\second{73.2\%} 
& 93.5\%/93.8\% 
& \second{67.0\%}/70.3\% 
& 40.4\%/43.4\%
& 65.8\%/70.9\% \\

OptiMUS   (2025)
& \best{86.9\%}/93.9\% 
& 54.0\%/60.0\% 
& 63.0\%/72.1\% 
& 87.8\%/89.0\% 
& 52.2\%/70.4\% 
& 38.5\%/41.3\%
& 63.7\%/71.1\% \\

OptiTree   (2025)
& 84.1\%/92.5\% 
& \second{64.0\%}/\second{71.0\%} 
& 62.5\%/71.9\% 
& 90.7\%/91.2\% 
& 65.4\%/\second{72.8\%} 
& 31.3\%/33.7\%
& 66.3\%/72.2\% \\

Lean-LLM-OPT   (2026)
& \second{86.5\%}/\second{94.9\%} 
& 55.0\%/63.0\% 
& 61.8\%/70.3\% 
& \second{93.8\%}/\second{94.2\%} 
& 57.6\%/70.0\% 
& \second{44.5\%}/\second{45.2\%}
& \second{66.5\%}/\second{72.9\%} \\

\midrule
\rowcolor{gray!20}
OptGraph (Ours)
& \best{86.9\%}/\best{95.3\%} 
& \best{71.0\%}/\best{81.0\%} 
& \best{66.0\%}/\best{78.2\%} 
& \best{93.9\%}/\best{94.4\%} 
& \best{76.4\%}/\best{89.2\%} 
& \best{58.4\%}/\best{60.8\%}
& \best{75.4\%}/\best{83.2\%} \\
\bottomrule
\end{tabularx}
\end{table*}

\subsubsection{Node \& Edge Update}

For each processed problem, OptGraph updates the current graph as
\begin{equation}
    \mathcal{G}^{(q)}
    =
    \mathcal{G}^{(q-1)}
    \oplus
    \Delta \mathcal{G}_q,
\label{eq:graph_update}
\end{equation}
where $\oplus$ denotes the insertion or merging of candidate nodes and edges. A candidate node is inserted when it represents previously unseen reusable knowledge. Otherwise, it is merged with a semantically and functionally similar node to reduce redundant graph growth. The corresponding relations are then added or updated to preserve the connections among problem types, modeling variants, formulations, implementations, errors, and repair knowledge.

\subsubsection{Final GraphRAG}

After all problems have been processed, OptGraph obtains the updated graph $\mathcal{G}^{(Q)}$.
The resulting GraphRAG preserves newly observed modeling patterns, including \texttt{Problem Type}, \texttt{Math Model}, \texttt{Typical Error}, and other nodes, as external knowledge.
This allows subsequent tasks to reuse experience accumulated from previous execution and validation processes without modifying the parameters of the underlying large language models.
In our implementation, the graph expands from $15$ initial seed patterns to $48$ final patterns after GraphRAG updating.

Overall, the procedure of OptGraph is summarized in Algorithm \ref{alg:optgraph_modeling}.
Different from conventional retrieval methods that return isolated examples or text snippets, OptGraph maintains a unified graph-structured modeling memory and updates it with execution traces and validation feedback.
It further updates this structured knowledge using execution traces and validation feedback, enabling non-parametric knowledge accumulation across tasks without modifying the parameters of the underlying LLMs.

\begin{figure*}[t]
    \centering

    \subfigure[ ]{
        \begin{minipage}{0.315\textwidth}
            \centering
            \includegraphics[width=\linewidth]{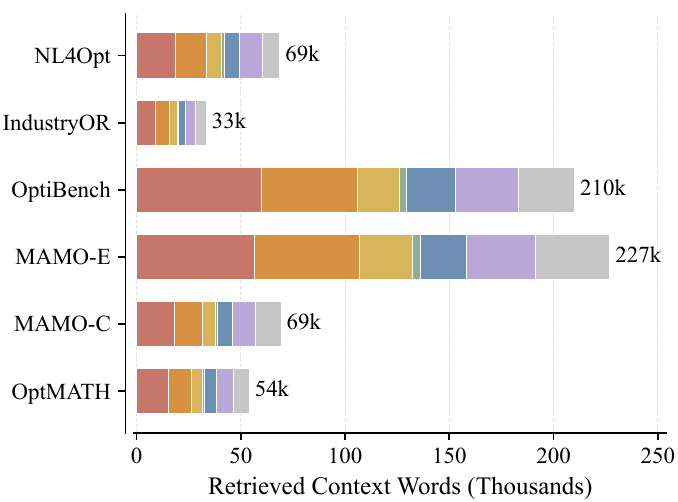}
            \vspace{-2mm}
            \label{fig:retrieval_volume}
        \end{minipage}
    }
    \hfill
    \subfigure[ ]{
        \begin{minipage}{0.315\textwidth}
            \centering
            \includegraphics[width=\linewidth]{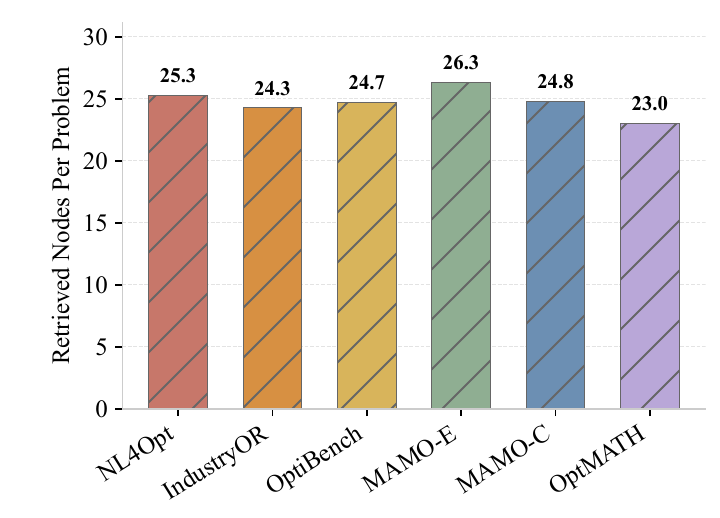}
            \vspace{-2mm}
            \label{fig:retrieval_nodes_per_problem}
        \end{minipage}
    }
    \hfill
    \subfigure[ ]{
        \begin{minipage}{0.315\textwidth}
            \centering
            \includegraphics[width=\linewidth]{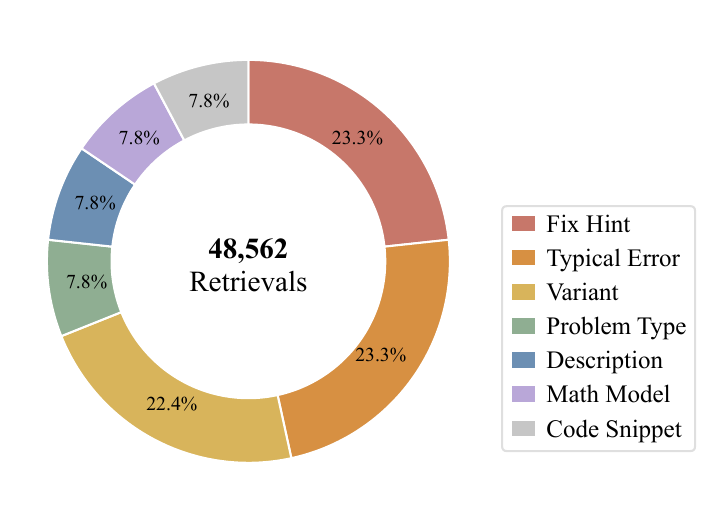}
            \vspace{-2mm}
            \label{fig:retrieval_node_distribution}
        \end{minipage}
    }

%    \vspace{-2mm}
    \caption{GraphRAG retrieval analysis across six datasets, where (a) retrieved context volume, (b) retrieved nodes per problem, (c) node-type distribution.}
    \label{fig:graphrag_retrieval_analysis}
    \vspace{-2mm}
\end{figure*}

\section{Experiments}
\label{sec:experiments}

This section evaluates the superiority of our OptGraph by comparing with the state-of-the-art methods such as CoE\footnote{\url{https://github.com/xzymustbexzy/Chain-of-Experts}}, OptiMUS\footnote{\url{https://github.com/teshnizi/OptiMUS}}, OptiTree\footnote{\url{https://github.com/MIRALab-USTC/OptiTree/tree/main}}, and Lean-LLM-OPT\footnote{\url{https://github.com/CoraLiang01/lean-llm-opt}}. Moreover, four LLMs including GPT-5.1, Gemini-3.1-Pro,  DeepSeek-V3.2, and Qwen3.6-Plus are added as baselines.

Subsection \ref{Experimental_Setups} introduces the experimental setups, Subsection \ref{Main_Results} gives the main results, Subsection \ref{Ablation_Studies} shows the ablation studies, 
Subsection \ref{Discussions} discusses the plug-and-play analysis, limitations of naive RAG, effects of validation rounds, and computational cost.

\begin{table*}[t]
\centering
\caption{Ablation studies (exact accuracy/5\% tolerance accuracy) of OptGraph on selected datasets, where \\the best results are marked in bold.}
\label{tab:ablation}
\renewcommand{\arraystretch}{1.2}
\footnotesize
\setlength{\tabcolsep}{5pt}
\begin{tabular}{ccccccccc}
\toprule
~~GraphRAG~~
& ~~Validation~~
& NL4Opt
& IndustryOR
& OptiBench
& MAMO-E
& MAMO-C
& OptMATH
& Average~~ \\
\midrule

$\times$
& $\times$
& 67.8\%/92.5\%
& 47.0\%/52.0\%
& 60.8\%/67.6\%
& 85.8\%/89.1\%
& 52.2\%/63.5\%
& 36.1\%/41.0\%
& 58.3\%/67.6\% \\

$\times$
& $\checkmark$
& 83.6\%/93.9\%
& 52.0\%/57.0\%
& 63.0\%/69.3\%
& 91.7\%/92.8\%
& 54.2\%/67.9\%
& 42.2\%/44.6\%
& 64.5\%/70.9\% \\

$\checkmark$
& $\times$
& 84.6\%/94.4\%
& 66.0\%/74.0\%
& 64.5\%/74.5\%
& 91.9\%/93.0\%
& 68.0\%/79.8\%
& 51.8\%/54.8\%
& ~~71.1\%/78.4\%~~ \\

\rowcolor{gray!20}
$\checkmark$
& $\checkmark$
& \textbf{86.9\%/95.3\%}
& \textbf{71.0\%/81.0\%}
& \textbf{66.0\%/78.2\%}
& \textbf{93.9\%/94.4\%}
& \textbf{76.4\%/89.2\%}
& \textbf{58.4\%/60.8\%}
& \textbf{75.4\%/83.2\%} \\
\bottomrule
\end{tabular}%
\end{table*}

\subsection{Experimental Setups}\label{Experimental_Setups}

\subsubsection{Datasets}

In the experiments, six representative optimization datasets with different sources, scales, and modeling difficulties are chosen. These datasets cover reusable OR modeling patterns such as objective interpretation, proportion and difference constraints, allocation, scheduling, routing, flow, blending, and other specialized structures. The details of these selected datasets are described as follows:
\begin{itemize}
\item NL4Opt\footnote{\url{https://github.com/NL4Opt/NL4Opt-competition}} contains 214 linear programming problems from the NeurIPS 2022 NL4Opt competition.
\item IndustryOR\footnote{\url{https://huggingface.co/datasets/CardinalOperations/IndustryOR}} has 100 real-world industrial OR problems, covering five categories and three difficulty levels.
\item OptiBench\footnote{\url{https://github.com/yangzhch6/ReSocratic}} includes 605 problems with diverse mathematical structures, such as linear, integer, nonlinear, and tabular optimization.
\item MAMO Easy (MAMO-E)\footnote{\url{https://github.com/FreedomIntelligence/Mamo}} collects 642 relatively simple linear programming problems with short descriptions and basic constraint structures.
\item MAMO Complex (MAMO-C)\footnote{\url{https://github.com/FreedomIntelligence/Mamo}} collects 203 complex problems with longer descriptions and richer constraints.
\item OptMATH\footnote{\url{https://github.com/optsuite/OptMATH}} contains  166 challenging modeling problems from mathematical programming instances.
\end{itemize}

\subsubsection{Evaluation Metrics}

In this paper, two popular metrics are considered, i.e., exact accuracy and 5\% tolerance accuracy. Specifically, exact accuracy follows the definition used in CoE \cite{xiao2024chain}, where a prediction is counted as correct only when the extracted answer exactly matches the ground-truth answer. For each problem, the generated code is executed, and the extracted answer $\hat{y}$ is compared with the ground-truth answer $y$. It is defined as
\begin{equation}
    \mathrm{Acc}_{\mathrm{exact}}
    =
    \frac{1}{N}
    \sum_{i=1}^{N}
    \mathbb{I}(\hat{y}_i = y_i),
\end{equation}
where $N$ is the number of problems
Besides, 5\% tolerance accuracy is adopted from OptiTree \cite{liu2025optitree}, which is defined as
\begin{equation}
    \mathrm{Acc}_{5\%}
    =
    \frac{1}{N}
    \sum_{i=1}^{N}
    \mathbb{I}
    \left(
    \frac{|\hat{y}_i-y_i|}{|y_i|+\eta}
    \leq 0.05
    \right),
\end{equation}
where  $\eta$ avoids division by zero. Obviously, higher precision accuracy and 5\% tolerance accuracy indicate better performance of  automated optimization.

\subsubsection{Implementation Details}
For compared methods, their released implementations are used directly, with only dataset interfaces, execution scripts, and answer extraction modules adapted when necessary. For OptGraph, the graph-structured knowledge base is initialized from reusable optimization resources and then used throughout the inference workflow. During inference, connected graph knowledge is retrieved to guide modeling, validation, repair, and graph updating.

We would like to point out that GPT-5.1 is utilized for Step 1, Step 2, and Step 4, covering GraphRAG construction, knowledge modeling, and GraphRAG update. Moreover,  Gemini-3.1-Pro is employed for Step 3, which serves as the verification mechanism. This setting  leverages a powerful model for graph construction, formulation generation, code generation, and knowledge updating, while adopting an independent model for validation to reduce same-model self-verification bias.

\begin{figure*}[t]
    \centering

    \subfigure[ ]{
        \begin{minipage}{0.315\textwidth}
            \centering
            \includegraphics[width=\linewidth]{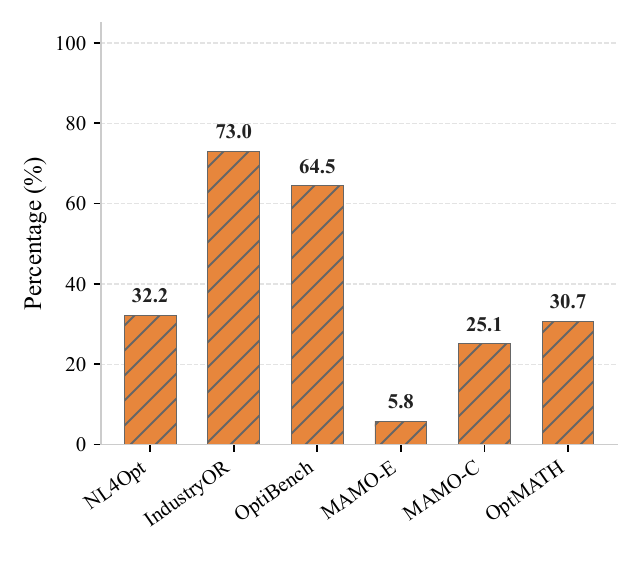}
            \vspace{-5mm}
            \label{fig:validation_trigger}
        \end{minipage}
    }
    \hfill
    \subfigure[ ]{
        \begin{minipage}{0.315\textwidth}
            \centering
            \includegraphics[width=\linewidth]{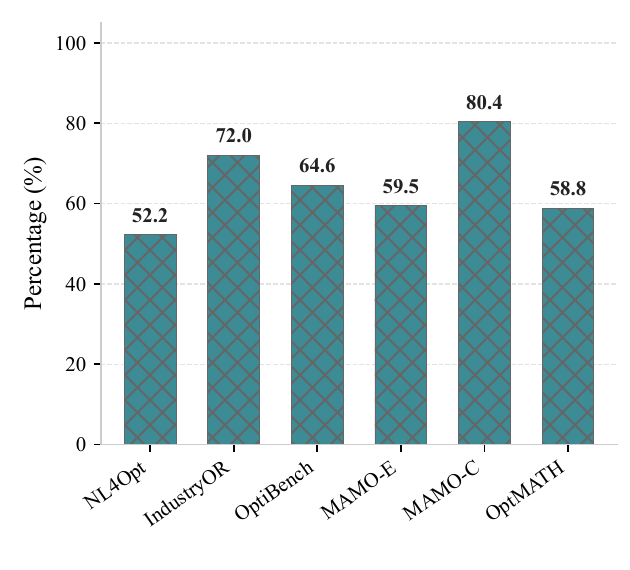}
            \vspace{-5mm}
            \label{fig:repair_success_exact}
        \end{minipage}
    }
    \hfill
    \subfigure[ ]{
        \begin{minipage}{0.315\textwidth}
            \centering
            \includegraphics[width=\linewidth]{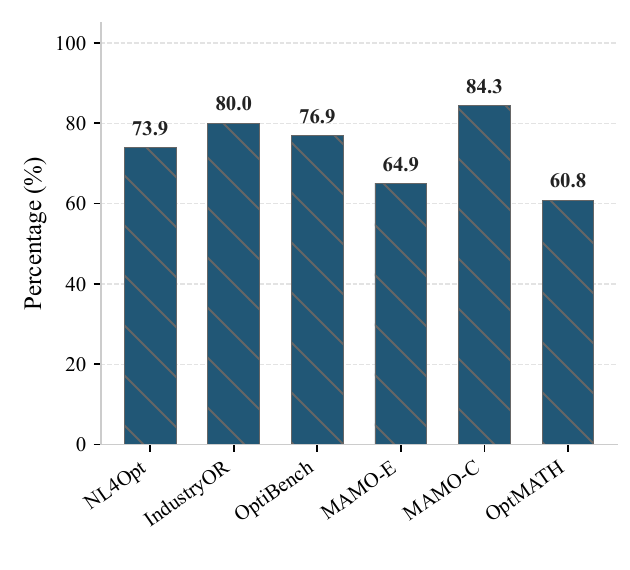}
            \vspace{-5mm}
            \label{fig:repair_success_5pct}
        \end{minipage}
    }

%   \vspace{-2mm}
    \caption{Effects of the validation mechanism across six datasets, where (a) validation trigger rate, (b) exact repair success rate, (c) 5\% repair success rate.}
    \label{fig:mechanism_analysis}
    \vspace{-2mm}
\end{figure*}

\begin{figure*}[t]
    \centering

    \subfigure[]{
        \begin{minipage}{0.315\textwidth}
            \centering
            \includegraphics[width=\linewidth]
            {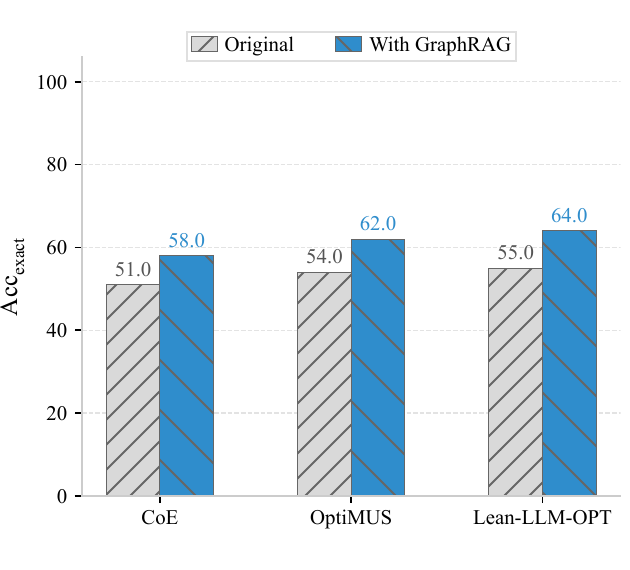}
            \vspace*{-4mm}
            \label{fig:plugin_industryor}
        \end{minipage}
    }
    \hfill
    \subfigure[]{
        \begin{minipage}{0.315\textwidth}
            \centering
            \includegraphics[width=\linewidth]
           {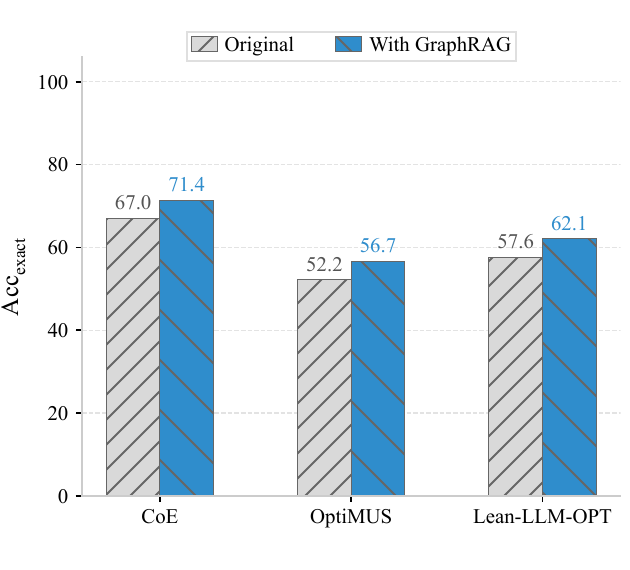}
            \vspace*{-4mm}
            \label{fig:plugin_mamo_complex}
        \end{minipage}
    }
    \hfill
    \subfigure[]{
        \begin{minipage}{0.315\textwidth}
            \centering
            \includegraphics[width=\linewidth]
        {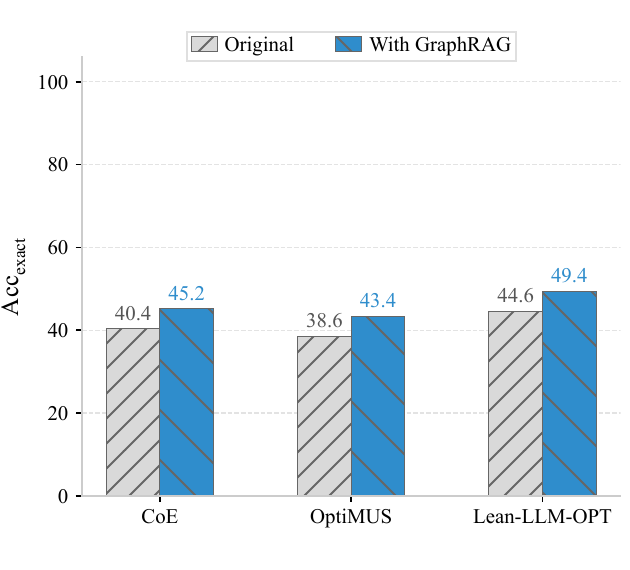}
            \vspace*{-4mm}
            \label{fig:plugin_optmath}
        \end{minipage}
    }

%    \vspace{-2mm}
    \caption{Plug-and-play analysis of GraphRAG on other LLM-based optimization  frameworks across  (a) IndustryOR, (b) MAMO-C,  (c) OptMATH.}
    \label{fig:plugin_effect}
    \vspace{-2mm}
\end{figure*}

\subsection{Main Results}\label{Main_Results}

Table \ref{tab3} reports the overall accuracy comparison on six datasets, with the best and second-best results highlighted in red and blue, respectively. In addition, the average results for each method are provided. It is evident that OptGraph achieves the highest exact accuracy and 5\% tolerance accuracy on all datasets, outperforming  CoE, OptiMUS, OptiTree, and Lean-LLM-OPT. Notably, OptGraph demonstrates a more pronounced advantage on complex datasets such as IndustryOR, MAMO-C, and OptMATH, where problems involve richer domain semantics, more complex constraints, and more challenging mathematical structures. For example, on OptMATH, the exact accuracy is improved by 13.9\% and 5\% tolerance accuracy by 15.6\% compared to the second-best method. This indicates that OptGraph can effectively handle complex optimization tasks through the introduction of GraphRAG and verification mechanisms. It is also found that, while standalone LLMs like GPT-5.1 and Gemini-3.1-Pro exhibit competitive optimization ability, but their performance remains clearly lower than agentic and graph-guided methods on average. This suggests that general-purpose LLMs alone are insufficient without structured modeling knowledge and validation support.

Furthermore, Fig. \ref{fig:graphrag_retrieval_analysis} illustrates the retrieval behavior behind these improvements. As demonstrated in Fig. \ref{fig:retrieval_volume}, the retrieved context volume varies across datasets and node types. Larger and more diverse datasets, such as OptiBench and MAMO-E, require more retrieved context, while IndustryOR and OptMATH use relatively compact but still structured retrieval. Fig. \ref{fig:retrieval_nodes_per_problem} shows that the average number of retrieved nodes per problem remains stable across datasets, suggesting that OptGraph does not simply increase retrieval size for harder datasets, but retrieves a compact graph neighborhood for each problem. Fig. \ref{fig:retrieval_node_distribution} reveals that the retrieved knowledge is not dominated by a single node type. Instead, OptGraph retrieves a balanced mixture of semantic, formulation, implementation, and repair-related knowledge.

Overall, these results demonstrate that our proposed OptGraph enhances automated optimization performance by retrieving knowledge from association graphs, thereby assisting LLMs in handling complex constraints and performing model code repair across IndustryOR, MAMO-C, and OptMATH.

\subsection{Ablation Studies}\label{Ablation_Studies}

Table \ref{tab:ablation} conducts the ablation experiments of OptGraph. 
It can be concluded that the full model achieves the best performance, confirming that graph-guided retrieval and validation-based repair work together to improve reliability. GraphRAG provides structured optimization knowledge for formulation and program generation, while the validation mechanism checks generated results and uses error-aware graph context to repair modeling, execution, and answer-extraction errors. Validation further checks the generated results and uses error-aware graph knowledge to repair formulation errors, execution failures, and answer-extraction mistakes.  
Removing either component leads to clear performance degradation, while removing both produces the weakest results across all datasets. 
Notably, on MAMO-C and OptMATH, removing both components results in a performance decline of over 20\% in exact accuracy and 5\% tolerance accuracy.

In addition, Fig. \ref{fig:mechanism_analysis} examines the effect of the validation mechanism. 
Fig. \ref{fig:validation_trigger}  displays the proportion of cases that trigger the validation mechanism. Obviously, IndustryOR and OptiBench have relatively high trigger rates, indicating that real-world and structurally diverse problems require more frequent checking and correction. In contrast, MAMO-E has a much lower trigger rate, which is consistent with its relatively simple and regular problem structure. Moreover, Fig. \ref{fig:repair_success_exact} and \ref{fig:repair_success_5pct} report the post-trigger repair success rates under exact and 5\% tolerance evaluation. The repair module achieves strong recovery on most datasets, especially under the 5\% tolerance criterion. These results confirm that validation is not merely an execution filter, but an active repair mechanism that works together with GraphRAG to correct formulation-level and code-level errors.

\subsection{Discussions}\label{Discussions}

This section discusses the plug-and-play analysis, limitations of naive RAG, effects of validation rounds, and computational cost on three complex datasets, i.e., IndustryOR, MAMO-C, and OptMATH.
 
\subsubsection{Plug-and-play Analysis}

This part evaluates whether GraphRAG can serve as a plug-and-play module for other LLM-based automated optimization frameworks. As shown in Fig. \ref{fig:plugin_effect}, the same graph-augmented retrieval mechanism is integrated into CoE, OptiMUS, and Lean-LLM-OPT. It can be found that GraphRAG consistently improves all three frameworks on these representative datasets. For example, on IndustryOR, the exact accuracy of CoE, OptiMUS, and Lean-LLM-OPT increases from 51.0\%, 54.0\%, and 55.0\% to 58.0\%, 62.0\%, and 64.0\%, respectively. Similar improvements are also observed on MAMO-C and OptMATH. These results convince us to believe that  graph-structured modeling knowledge can be incorporated into different prompt-based or agentic optimization modeling systems.

\begin{figure}[t]
    \centering
    \includegraphics[width=0.85\linewidth]{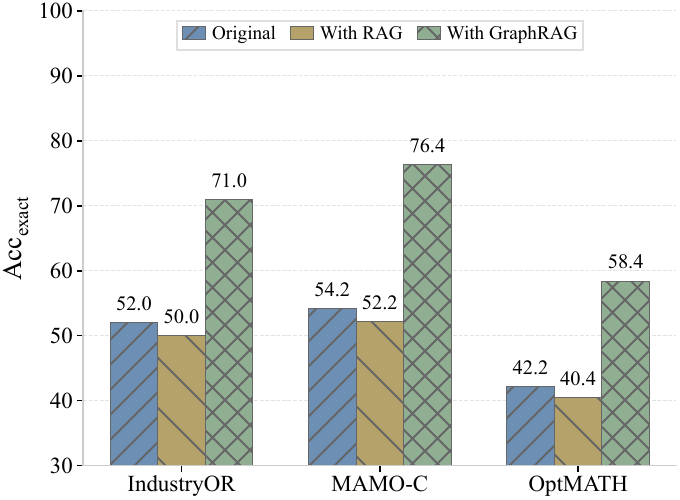}
    \caption{Effects of naive RAG and GraphRAG on the original workflows across IndustryOR, MAMO-C, and OptMATH.}
    \label{fig:vector_rag_comparison_3datasets}
\end{figure}

\begin{figure}[t]
    \centering
    \includegraphics[width=0.85\linewidth]{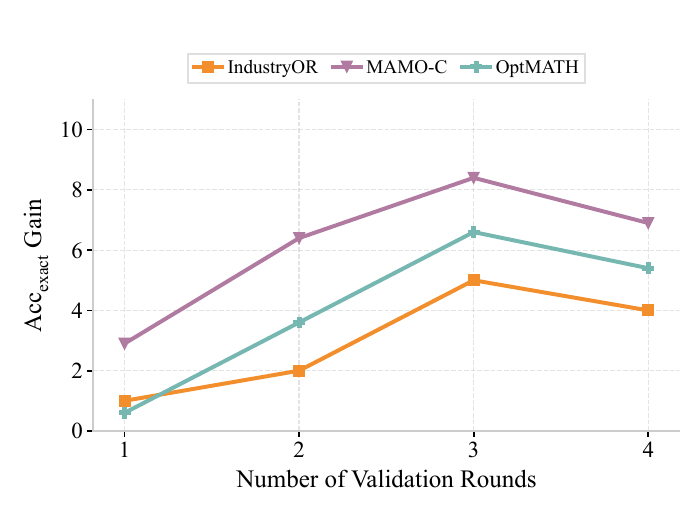}
    \caption{Accuracy gains over the zero-round setting when varying the number of validation rounds across IndustryOR, MAMO-C, and OptMATH.}
    \label{fig:validation_rounds_gain}
\end{figure}

\subsubsection{Limitations of Naive RAG}

Fig. \ref{fig:vector_rag_comparison_3datasets} compares the original workflows with naive RAG and GraphRAG on the same  datasets. The results show that simply adding retrieved examples does not necessarily improve optimization performance. On all three datasets, naive RAG performs worse than the original workflows, decreasing exact accuracy from 52.0\% to 50.0\% on IndustryOR, from 54.2\% to 52.2\% on MAMO-C, and from 42.2\% to 40.4\% on OptMATH. In contrast, our OptGraph with GraphRAG achieves 71.0\%, 76.4\%, and 58.4\% on these datasets, respectively. This suggests that semantic similarity alone may retrieve examples that are textually similar but not aligned with the underlying logic. Therefore, retrieval for automated optimization should not only return similar text snippets, but also preserve structured relations.

\subsubsection{Effects of Validation Rounds}

Fig. \ref{fig:validation_rounds_gain} analyzes the effect of the maximum number of validation rounds. 
The reported values denote the exact-accuracy gain over the zero-round setting.  Increasing the validation budget generally improves performance on all three datasets, showing that iterative validation and repair are useful for complex optimization modeling tasks. 
The gains peak at 3 validation rounds, where IndustryOR, MAMO-C, and OptMATH obtain improvements of 5.0\%, 8.4\%, and 6.6\%, respectively.  When the budget is increased to 4 rounds, the gains decrease on all three datasets, indicating that excessive repair may introduce unnecessary changes or over-correction. Therefore, the maximum number of validation rounds is set to 3 by default, which provides a good balance between accuracy improvement and computational cost.

\subsubsection{Computational Cost}

Table \ref{tab:efficiency_accuracy} lists the computational costs and average accuracy of the different frameworks. It can be seen that our proposed OptGraph consumes more tokens than lightweight methods such as CoE and Lean-LLM-OPT because it retrieves graph-structured modeling knowledge and performs validation-based repair. Meanwhile, the average runtime of OptGraph is also higher than that of OptiTree and Lean-LLM-OPT.
However, among all the compared methods, OptGraph achieves the highest average exact accuracy and 5\% tolerance accuracy. It is noted that compared with OptiMUS, our OptGraph requires fewer average seconds per problem while achieving substantially higher accuracy. These results demonstrate that OptGraph trades off additional retrieval and verification costs for greater formulation robustness and solution reliability, making it more suitable for accuracy-oriented automated optimization.

\begin{table}[!t]
\centering
\caption{Efficiency and average accuracy comparisons, where token usage is estimated from visible logs, and the best\\ results are marked in bold.}
\label{tab:efficiency_accuracy}
\renewcommand{\arraystretch}{1.15}
\footnotesize
\setlength{\tabcolsep}{3pt}
\begin{tabular}{lcccc}
\toprule
Method
& ~~Tokens (k)~~
& ~~Time (s)~~
& ~~~$\mathrm{Acc}_{\mathrm{exact}}$~~~
& ~~$\mathrm{Acc}_{5\%}$~~ \\
\midrule

CoE
& \textbf{1.3}
& 79.1
& 65.8
& 70.9 \\

OptiMUS
& 15.9
& 111.1
& 63.7
& 71.1 \\

OptiTree
& 4.6
& \textbf{14.3}
& 66.3
& 72.2 \\

Lean-LLM-OPT
& 1.5
& 17.1
& 66.5
& 72.9 \\

\rowcolor{gray!20}
OptGraph
& 9.7
& 93.9
& \textbf{75.4}
& \textbf{83.2} \\
\bottomrule
\end{tabular}%
\end{table}

\section{Conclusion}
\label{sec:conclusion}

This paper proposed a novel graph-augmented agentic workflow for automated optimization. Instead of relying on fixed prompting strategies or isolated retrieved examples, OptGraph leverages reusable optimization knowledge as a structured graph that connects \texttt{Problem Type}, \texttt{Variant}, \texttt{Description}, \texttt{Math Model}, \texttt{Code Snippet}, \texttt{Typical Error}, and \texttt{Fix Hint}. Moreover, the graph is adaptively updated based on the execution trajectories and validation feedback. Numerical studies show that OptGraph achieves competitive or even superior accuracy. These findings suggest that structured and reusable modeling knowledge is a promising direction for building more robust LLM-based evolutionary optimization systems.

In the future, more efficient graph-augmented frameworks can be developed to reduce token consumption and shorten inference time while preserving modeling accuracy \cite{xiao2025survey}. Furthermore, in industrial-scale scenarios involving vast numbers of variables and multi-stage stochastic structures, existing LLM-based methods fall short of expectations regarding accuracy and stability \cite{da2026large}.
Finally, the generated graph-structured trajectories facilitate the analysis of how retrieved  knowledge influences problem formulation decisions and can be further integrated with symbolic verification techniques \cite{wang2025translating}.

\bibliographystyle{IEEEtran}
\bibliography{mybib}

\end{document}